\documentclass[twoside, 12pt]{report}
\usepackage[english]{babel}
\usepackage{epsfig}
\usepackage{psfig}
\usepackage{amssymb,amsmath,amsfonts,soul,amsthm,enumerate}
\usepackage{graphicx}
\usepackage[cp1251]{inputenc}
\usepackage[T2A]{fontenc}
\usepackage{mathrsfs}
\usepackage{cite} 

\usepackage{tikz}

\tikzset{square arrow/.style={
		to path={-- ++(0,-.45)  -| (\tikztotarget) \tikztonodes},below,pos=.25}}

\def\@evenfoot{}
\def\@oddfoot{}

\usepackage{authblk}

\begin{document}

\def\@evenhead{\vbox{\hbox to \textwidth{\thepage\leftmark}\strut\newline\hrule}}

\def\@oddhead{\raisebox{0pt}[\headheight][0pt]{%
\vbox{\hbox to \textwidth{\rightmark\thepage\strut}\hrule}}}

\def\bibname{\vspace*{-30mm}{\centerline{\normalsize References}}}

\newpage
\normalsize
\thispagestyle{empty}

\vskip 5 mm

\centerline{\bf FUNCTIONAL-ANALYTICAL RECONSTRUCTION}
\centerline{\bf OF HIGH CONTRAST INHOMOGENEITIES}

\vskip 0.3cm
\centerline{\bf A.S. Shurup$^{ \textbf{1),2)}}$}

\vskip 0.3cm

\begin{center}
\begin{minipage}{118mm}
	\small
	\textit{$^{1)}$~M.V.~Lomonosov Moscow State University, Faculty of Physics, Acoustics Department,
		Leninskie Gory, Moscow 119991, Russia}.
		
	\textit{$^{2)}$~Sсhmidt Institute of Physics of the Earth
		of the Russian Academy of Sciences, 
		B. Gruzinskaya str., 10, build. 1, Moscow 123242, Russia.}

	\centerline{ \textit{e-mail: shurup@physics.msu.ru}}

\end{minipage}
\end{center}

\vskip 0.3cm
\noindent{\bf Abstract}
{\small In practice of acoustic tomography, for example, in medical applications and ocean tomography, the relative deviation of sound speed from its background value usually does not exceed 10-30\%. At the same time, in electromagnetic applications, the equivalent contrasts can be noticeably higher than 60\%. Since the inverse electromagnetic problem can be reduced in some approximation to Helmholtz equation, a formal comparison of reconstruction results obtained for different "acoustic" contrast and corresponding "dielectric" contrast is possible. In this work examples of such reconstructions are presented, which were obtained by using the functional-analytical algorithm described in works of R.G.~Novikov. Previously, the advantages of this algorithm for solving practical problems of acoustic tomography were demonstrated. Results obtained in the present work show that functional-analytical algorithm can also be applied to reconstructing inhomogeneities with high "dielectric" contrast. Moreover, the functional algorithm also perfectly reconstructs very small "dielectric" contrast, recovering of which can be difficult for other approaches due to weak backscattering.
\vskip 0.2cm
\noindent {\bf Key words:} inverse acoustic scattering, numerical modeling
\vskip 0.2cm
\noindent {\bf AMS Mathematics Subject Classification:}  35R30, 65N21}
\vskip 0.3cm

\setcounter{figure}{0}

\renewcommand{\thesection}{\large 1}

\section{\large Introduction}

Acoustic tomography allows to reconstruct internal structure of almost any natural environment with  appropriate choice of frequency range of probing signal. Among actual areas of acoustic tomography application are medicin \cite{label2,label1}, ocean tomography  \cite{label3,label4}, as well as tomographic methods for studying the Earth \cite{label5,label6}. From mathematical point of view, acoustic tomography is a special case of a more general class of inverse scattering problems. Mathematically rigorous functional-analytical methods for solving inverse problems are known \cite{label7,label8,label9,label9_1,label10,label11,label12,label13}, which were initially developed for quantum mechanical applications. One of such algorithm, proposed in works of R.G. Novikov \cite{label14,label15}, turned out to be very promising for solving practical problems of acoustic tomography that was demonstrated in detailed numerical investigations \cite{label16,label17,label18,label19,label20,label21}. Some advantages of this algorithm include the following: multiple scattering processes are taken into account that allows reconstruction of scatterers beyond the Born approximation without iterations and additional regularization \cite{label16, label18}; possibility for the joint reconstruction of scalar-vector inhomogeneities, which includes space distributions of sound speed, density, frequency dependent attenuation and vector field of flows \cite{label19, label22}; generalization of the algorithm to multifrequency regime of sounding, which provides the high interference resistance of this algorithm, which is necessary for practical applications \cite{label17}; scatterer functions can be found in different space points independently \cite{label16} that gives wide opportunities for applications of parallel computation methods; multi-channel scattering effects can be taken into account during inverse problem solution for reconstruction of 3D inhomogeneity located in a layer, for example in ocean waveguide \cite{label20, label21}. The aforementioned advantages distinguish the considered functional-analytical algorithm from other methods of solving inverse problems of acoustic scattering. At the same time, the question about possibility of reconstructing inhomogeneities with high contrast using this approach was remaining open. The point is that in practice of acoustic tomography (in medical applications, ocean tomography) deviations of sound speed $c$ from background value $c_0$ do not exceed $ |c - c_0| / c_0 = | \Delta c | \big/ c_0 \le 0.1 \div 0.3 $ \cite{label23, label24, label25} that are relatively small contrast. At the same time, in electromagnetic applications \cite{label26, label27} the equivalent contrasts can be noticeably higher. In \cite{label26, label27} the inverse electromagnetic problem is reduced under some approximations to Helmholtz equation and effective method of reconstructing inhomogeneities with high "dielectric" contrast was developed and validated on experimental data. "Dielectric" contrast, or dielectric constant, is the ratio of permittivity inside an inhomogeneity to permittivity of background medium \cite{label27_1}. Thus, a formal application of discussed functional-analytical algorithm for reconstructing high "dielectric" contrasts (significantly exceeding the corresponding values $ | \Delta c | \big/ c_0 \approx 0.3$ ) is possible. In this work examples of such reconstructions are presented both for the case $c \gg c_0$ and  $c \ll c_0$, demonstrating the abilities of functional-analytical reconstruction in these cases.

\renewcommand{\thesection}{\large 2}
\section{\large Acoustic and dielectric contrasts}

\begin{figure}[ht]
	\centerline{\epsfig{file=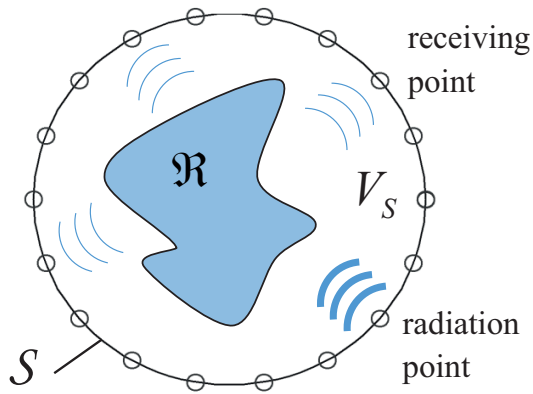, width = 9 cm}}
	\caption{Tomography area $V_S$ contains the scattering region  $\mathfrak{R}$; quasi-point sources and receivers are located on the boundary $S$.}
	\label{figure_1}
\end{figure}

The two-dimensional inverse problem is considered. 
It is assumed that on the boundary $ S $  of tomography area $ V_S $  there are transducers, which are equivalent to point ones (quasi-point transducers), emitting and receiving acoustic fields. 
Inside the region $ V_S $ there is an inhomogeneity (scatterer), which is nonzero only inside the scattering domain $\mathfrak{R}$, which lies entirely inside $ V_S $: $\mathfrak{R} \subset V_S$ (figure~\ref{figure_1}). 
For simplicity, those inhomogeneities will be considered, which are described by sound speed perturbations only. A more general case of joint reconstruction of sound speed, density, absorption coefficient and flows can also be implemented on the basis of relations given in \cite{label22}.

The spatial distribution of complex spectral amplitude of acoustic pressure $ p(\mathbf{r}) $ in the considered inhomogeneous medium, which is characterized by the sound speed  $c(\mathbf{r})$, constant density of a medium without flows and with negligible absorption, is described by the Helmholtz equation \cite{label28}:
\begin{equation}
	\label{Helm_eq}
	\nabla^2 p(\mathbf{r}) + \frac{\omega^2}{c^2(\mathbf{r})} \ p(\mathbf{r}) = 0,
\end{equation} 
\noindent
here $\mathbf{r}$ is a radius vector, $\omega$ is a circular frequency, time dependence is assumed in the form $\sim \text{exp}(- i \omega t)$. 

As it was shown in \cite{label26, label27}, the process of electromagnetic wave scattering on inhomogeneities of dielectric constant $\epsilon(\mathbf{r})$ can also be described by the Helmholtz equation, which has the form:
\begin{equation}
	\label{Helm_eq_dielect}
	\nabla^2 p(\mathbf{r}) + \frac{\omega^2}{c^2_{00}} \epsilon(\mathbf{r}) \ p(\mathbf{r}) = 0,
\end{equation} 
\noindent
here $c_{00}$ is the speed of light in background medium; it is assumed that the permittivity of the background medium is 1, i.e. the value $\epsilon(\mathbf{r})$ characterizes the "dielectric" contrast of the inhomogeneity in question. 

From (\ref{Helm_eq}), (\ref{Helm_eq_dielect}) a formal acoustic analogy for $\epsilon(\mathbf{r})$ follows 
\begin{equation}
	\label{contrasts}
	\epsilon(\mathbf{r}) = \frac{c^2_0}{c^2(\mathbf{r})} = \left( 1 + \frac{\Delta c(\mathbf{r})}{c_0} \right)^{-2},
\end{equation} 
\noindent
which connects "dielectric" $\epsilon(\mathbf{r})$ and "acoustic" $\Delta c(\mathbf{r}) \big/ c_0 \equiv ( c(\mathbf{r}) - c_0 ) \big/ c_0$ contrasts. For example, if sound speed perturbation $\Delta c(\mathbf{r}) \big/ c_0 $ = 0.3, which is quite large for acoustic tomography, then "dielectric" contrast is quite small $\epsilon(\mathbf{r}) \approx 0.59$. In \cite{label26, label27} a method for solving the inverse problem was proposed and tested on experimental data, which makes it possible to reconstruct high "dielectric" contrasts $\epsilon \gg 1$. It is important to investigate reconstructing similar contrasts by using algorithms, whose possibilities in acoustic inverse problems have already been demonstrated in the framework of numerical simulation.

Equation (\ref{Helm_eq}) can be rewritten in the form:
\begin{equation}
	\label{Helm_eq_LS}
	\nabla^2 p(\mathbf{r}) + k^2_0 \ p(\mathbf{r}) = v(\mathbf{r}) \  p(\mathbf{r}) ,
\end{equation} 
\noindent
here $k_0$ is the wavenumber of background medium and the scatterer function $v(\mathbf{r})$ has the form:
\begin{equation}
	\label{scatt_func}
	v(\mathbf{r}) = k^2_0 \left[ 1 - \epsilon(\mathbf{r}) \right].
\end{equation} 

Application of the Novikov’s functional-analytical algorithm \cite{label14,label15} allows to find an estimate $\hat{v}(\mathbf{r})$ of unknown scatterer function $v(\mathbf{r})$ in (\ref{Helm_eq_LS}) by using acoustic fields $p(\mathbf{r})$ radiated and received at the boundary of the tomographic area. Scatterer function (\ref{scatt_func}) will be presented below as a results of reconstruction.

\renewcommand{\thesection}{\large 3}
\section{\large Numerical modeling}

In the numerical simulation, a two-dimensional region $V_S$ of cylindrical shape, surrounded by 40 receiving-emitting transducers was considered. The radius of region was assumed to be equal to $R = 7 \lambda_0$, where $\lambda_{0}$ is the wavelength in the background environment, expressed in relative length sampling units (l.s.u.): $\lambda_0$~=~8~l.s.u. The center of reconstructed inhomogeneities was shifted relative to the center of tomography region along $ Ox $ axis by $\lambda_0 / 2$ to eliminate the symmetry of the considered problem. 

Since sources and receivers are assumed to be quasi-point transducers, acoustic fields in an arbitrary point $\mathbf{r}$ are the Green's functions  $G_0(\mathbf{r}, \mathbf{x})$ or $G(\mathbf{r}, \mathbf{x})$ in the absence or in the presence of a scatterer, respectively; here $\mathbf{x}$ means the radius vector of a source. In the considered two-dimensional case $G_0(\mathbf{r}, \mathbf{x}) = - (i/4) \, H^{(1)}_0(k_0 |\mathbf{r} - \mathbf{x}|)$, where $ H^{(1)}_0 $ is the Hankel function of the zero order and the first kind.

In this work, the reconstruction of sound speed inhomogeneities of cylindrical shape is considered. The scattering of a cylindrical wave on a cylindrical inhomogeneity has a rigorous analytical solution (see, for example, \cite{label31}), that allows to calculate scattering data $G(\mathbf{r}, \mathbf{x})$ with very high quality, in comparision with a discretized solution of the Lippmann-Schwinger equation. In the latter case, the question arises about the spatial discretization step of the integral equation, which must decrease as the scatterer strength increases in order to adequately describe the processes of multiple scatterings. When reconstructing inhomogeneities with high contrasts, this leads to very small step and, as a result, causes restrictions on computer memory, thereby limiting the possibilities of modeling. When analytical solution of the scattering problem is used for cylinders with high contrasts, such limitations do not arise. 

\begin{figure}[ht!]
	\centerline{\epsfig{file=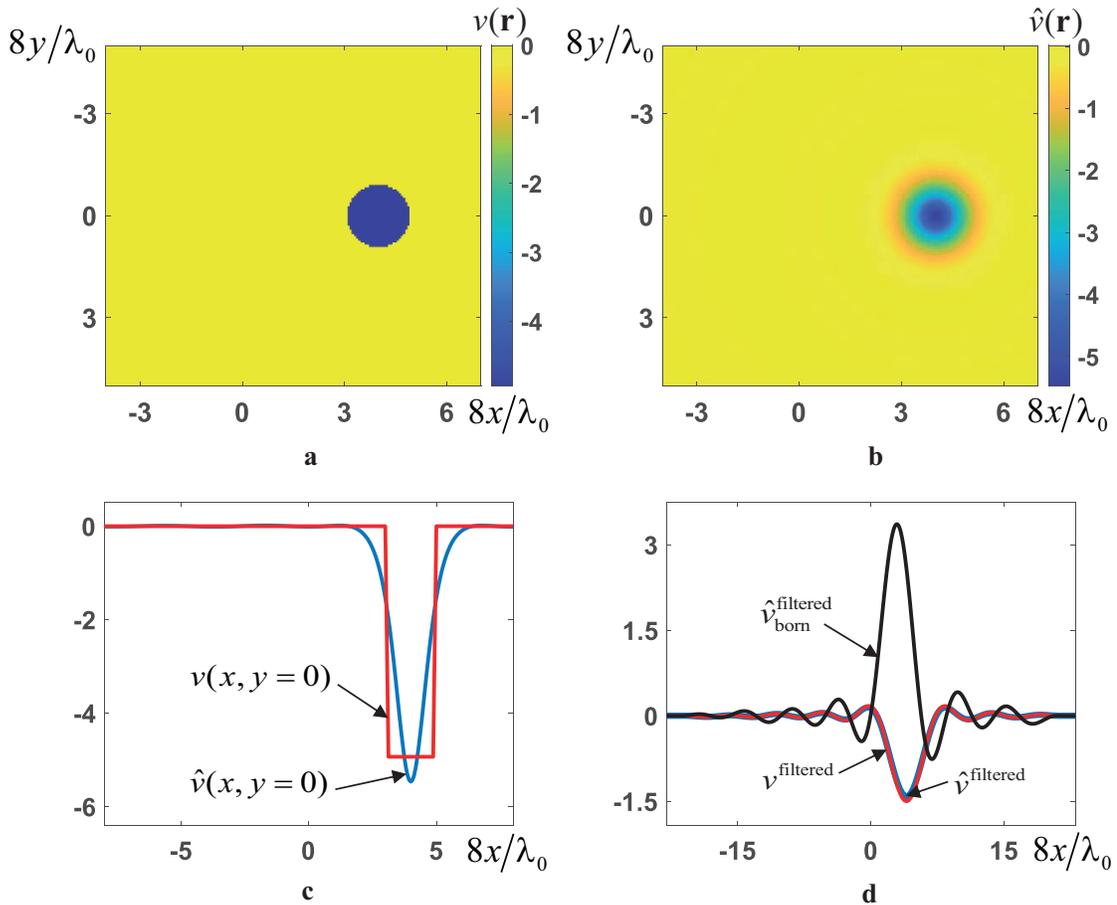, width = 15.5 cm}}
	\caption{The case $c \ll c_0 $, $\epsilon \gg 1$. General view of true scatterer $v(\mathbf{r})$ (a) of cylindrical shape with relative contrast of sound speed $\Delta c \big/ c_0 \approx - 0.67$ and corresponding "dielectric" contrast  $\epsilon = 9$, which provides the additional phase shift $\Delta \psi \approx 0.94 \pi$ and norm of scattering data $\| f(\phi, \phi^\prime) \| \approx 8.6 / (3 \pi)$, relative size of inhomogeneity is $ 2 R_0 / \lambda_0 \approx 0.25$;
		\protect\\
		- space distribution of reconstruction results $\hat{v}(\mathbf{r})$ obtained by using the functional-analytical algorithm (b); 
		\protect\\
		- central cross sections of true scatterer $v(x, y = 0)$ (c, red line) and of estimate $\hat{v}(x, y~=~0)$ (c, blue line); 
		\protect\\
		- central cross sections of functions which space spectrums are filtered in a cylinder with radius $2k_0$ (arguments of functions are omitted): filtered true scatterer $ v^{\text{filtered}}$ (d, red line), functional-analytical estimate $ \hat{v}^{\text{filtered}}$ (d, blue line) and Born approximation result $ \hat{v}^{\text{filtered}}_{\text{born}}$ (d, back line).}
	\label{figure_2}
\end{figure}

However, the following peculiarity should be taken into account. An inhomogeneity of a cylindrical shape has a wide spatial spectrum, which goes beyond boundaries of the circle with radius $2 k_0$ in a space of wave vectors. It is known \cite{label16} that in monochromatic two-dimensional inverse scattering problem it is possible to reconstruct spatial frequencies of a scatterer, localized mainly inside the circle $2 k_0$. Therefore, reconstruction results given below do not provide a cylinder with sharp boundaries, but the result of spatial filtering of this cylinder (figures~\ref{figure_2},~\ref{figure_3}). This feature should be taken into account when analyzing  reconstruction results presented below. 

To characterize the accuracy of estimates $\hat{v}(\mathbf{r})$, the relative root-mean-square (rms) reconstruction errors are calculated over the entire tomography region $V_S$: 
\[
\delta_v \equiv \sqrt{\int\limits_{V_S} \bigl| \hat{v}(\mathbf{r}) - v(\mathbf{r}) \bigr|^2 d \mathbf{r}} \bigg/ \sqrt{\int\limits_{V_S} \bigl| v(\mathbf{r}) \bigr|^2 d \mathbf{r}}.
\]

To describe the strength of considered scatterers, i.e. to estimate how strongly they distort the incident acoustic field, the values of additional phase shifts are calculated
\[
\Delta \psi = k_0 \int\limits_{l_\mathfrak{R}} \frac{\Delta c(\mathbf{r}) \big/ c_0}{1 + \Delta c(\mathbf{r}) \big/ c_0} d l_\mathbf{r} = k_0 \int\limits_{l_\mathfrak{R}} \bigl[ 1 - \sqrt{\epsilon(\mathbf{r})}  \bigr] d l_\mathbf{r},
\]

\noindent
which are caused due to a wave propagation along trajectory $l_\mathfrak{R}$ through scatterer with relative "acoustic" speed contrast $\Delta c(\mathbf{r}) \big/ c_0 $ or with "dielectric" contrast $\epsilon(\mathbf{r})$ (\ref{contrasts}); $d l_\mathbf{r}$ is the length of trajectory element in vicinity of point $\mathbf{r}$. The norm of scattering data is also calculated in the form
\[
\| f(\phi, \phi^\prime) \| \equiv \sqrt{\int\limits_0^{2 \pi} d \phi \int\limits_0^{2 \pi} d \phi^\prime \ \bigl| f(\phi, \phi^\prime) \bigr|^2 },
\]
\noindent
where $\phi$, $\phi^\prime$ are angular components of wave vectors of incident $\boldsymbol{\kappa} = \left\{ k_0, \phi \right\}$, and scattered $\boldsymbol{\ell} = \left\{ k_0, \phi^\prime \right\} $ plane waves; $f(\phi, \phi^\prime)$ is the classical scattering amplitude, characterizing scattered fields in far zone \cite{label18, label19}. The norm of scattering data $\| f(\phi, \phi^\prime) \|$, together with the additional phase shift $\Delta \psi$, characterizes the scatterer strength. Recalculation of initial scattering data $G(\mathbf{r}, \mathbf{x})$ to scattering amplitude $f(\phi, \phi^\prime)$ is implemented on the basis of relations presented in \cite{label12,label18,label29}; an alternative way of such recalculation can be found in \cite{label30}.

Figure~\ref{figure_2} shows reconstruction results obtained for scatterer (5) with contrasts $\epsilon~=~9$, $\Delta c \big/ c_0 \approx -0.67$ and with relative size $ 2 R_0 / \lambda_0 \approx 0.25$. For problems of ocean tomography, the common value of sound speed in background medium is $c_0 = 1500 $~m/s, then in the considered case sound speed inside inhomogeneity is $c = 500$~m/s, i.e. three times less than in a water under normal conditions. Such unusual case $c \ll c_0$ can occur in practice, for example, in shallow waters with a subsurface layer of bottom containing methane bubbles (gas-saturated bottom) \cite{label32}. The presence of methane can be caused, for example, by processes of gas release in areas of oil-and-gas basin on shelf of northern seas, thus, the detection of such sound velocity anomalies in bottom can be used for hydrocarbon deposits location. For the considered inhomogeneity, the additional phase shift is $\Delta \psi \approx 0.94$ and the norm of scattering amplitude is $\| f(\phi, \phi^\prime) \| \approx 8.6 / (3 \pi)$. This scatterer is strong enough and cannot be reconstructed with acceptable accuracy within the Born approximation. As it can be seen in figures~\ref{figure_2}a-\ref{figure_2}c, the location of inhomogeneity and its amplitude value are reconstructed by the functional-analytical algorithm with acceptable accuracy. But the form of obtained estimate $\hat{v}(\mathbf{r})$ is smoother than the original cylinder due to the space filtration, as it was noted above. This leads to overestimated value of residual $\delta_v \approx 0.47$. After spatial filtering in the circle $2k_0$, initial scatterer $v^{\text{filtered}}(\mathbf{r})$ and reconstructed function $ \hat{v}^{\text{filtered}}(\mathbf{r})$ become almost visually indistinguishable $v^{\text{filtered}}(\mathbf{r}) \simeq \hat{v}^{\text{filtered}}(\mathbf{r})$ (figure~\ref{figure_2}d); the value of discrepancy $\delta_v$ in this case decreases to 0.05. It should be noted that an attempt to reconstruct the same scatterer in the Born approximation (i.e., in the single scattering approximation) gives unsatisfactory results (figure~\ref{figure_2}d) – in particular, even after spatial filtering the discrepancy is equal to $\delta_v \approx 3.1$. During numerical modeling, inhomogeneities with other contrasts $\epsilon = 2, 5, 8$, corresponding to sound speed perturbations $\Delta c \big/ c_0 \approx -0.29, -0.55, -0.65$, have also been reconstructed; obtained results are close in quality to reconstructed estimates shown in figure~\ref{figure_2} and are not presented here to reduce the amount of graphic materials.

\begin{figure}[ht!]
	\centerline{\epsfig{file=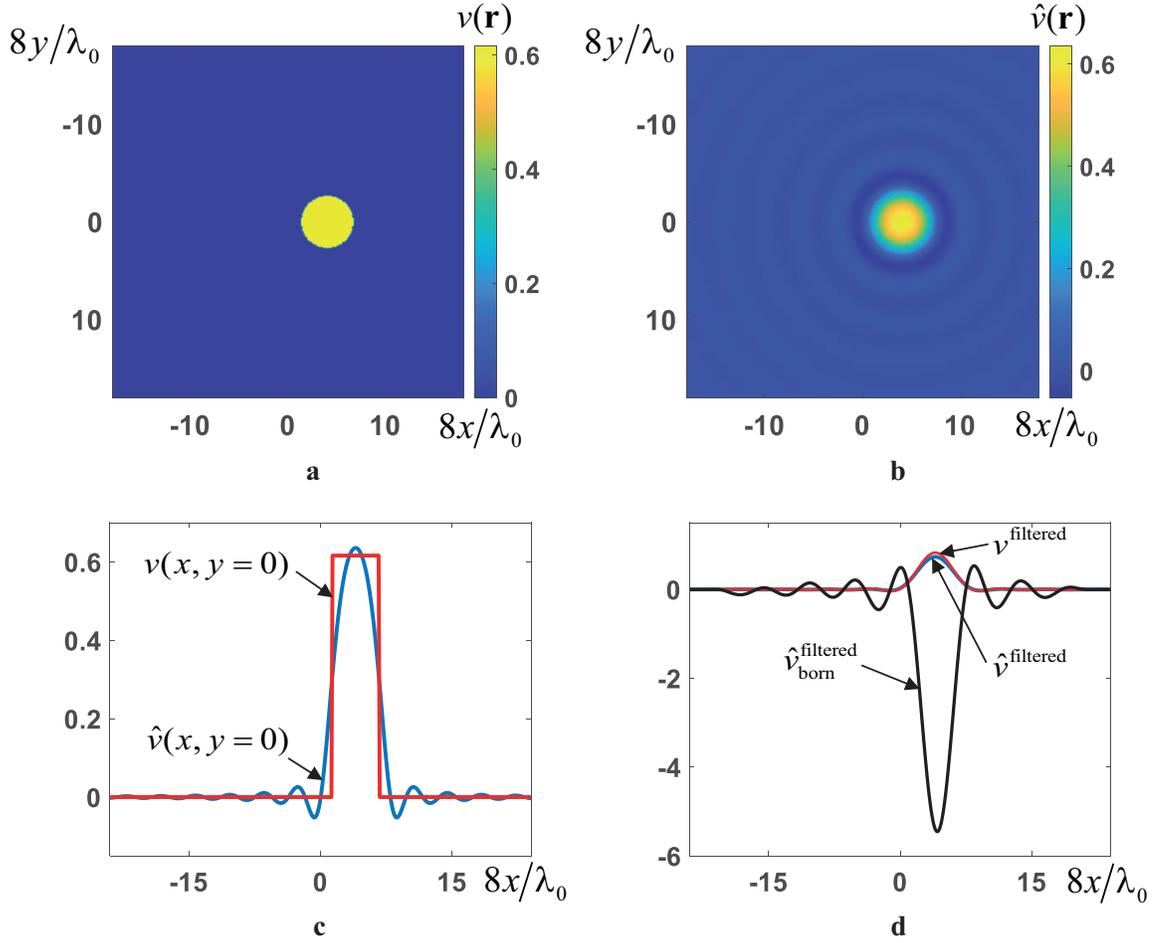, width = 15.5 cm}}
	\caption{The case $c \gg c_0 $, $\epsilon \ll 1$. Space distribution of true scatterer $v(\mathbf{r})$ (a) with relative contrast of sound speed  $\Delta c \big/ c_0 = 39$ and corresponding "dielectric" contrast  $\epsilon \approx 6 \cdot 10^{-4}$, which provides the additional phase shift $\Delta \psi \approx 1.3 \pi$ and norm of scattering data $\| f(\phi, \phi^\prime) \| \approx 7.7 / (3 \pi)$, relative size of scatterer is $ 2 R_0 / \lambda_0 \approx 0.67$;
		\protect\\
		- general view of functional-analytical estimate  $\hat{v}(\mathbf{r})$ (b); 
		\protect\\
		- central cross sections of true inhomogeneity $v(x, y = 0)$ (c, red line) and reconstruction result  $\hat{v}(x, y~=~0)$ (c, blue line); 
		\protect\\
		- central cross sections ($y = 0$) of space filtering results obtained in a cylinder with radius  $2k_0$ (arguments of functions are omitted): filtered true scatterer $ v^{\text{filtered}}$ (d, red line), functional-analytical estimate $ \hat{v}^{\text{filtered}}$ (d, blue line) and Born approximation result $ \hat{v}^{\text{filtered}}_{\text{born}}$ (d, back line).}
	\label{figure_3}
\end{figure}

Figure~\ref{figure_3} shows reconstruction results for the case $c \gg c_0$, which may correspond to detection of solid objects under water, for example, sunken objects or ships during archaeological research. To test the possibilities of considered algorithm, an enormously large (and thus non-physical) contrast was considered $\Delta c \big/ c_0 = 39$, which corresponds to very small values of "dielectric" contrast $\epsilon \approx 6 \cdot 10^{-4}$. For the considered case relative size of inhomogeneity is $ 2 R_0 / \lambda_0 \approx 0.67$ that gives additional phase shift $\Delta \psi \approx 1.3$ and norm of scattering data $\| f(\phi, \phi^\prime) \| \approx 7.7 / (3 \pi)$. This inhomogeneity is also sufficiently strong, as it was before for the case $c \ll c_0$. In figures~\ref{figure_3}a-\ref{figure_3}c one can see, that location of scatterer function and its amplitude are reconstructed by the functional-analytical algorithm with acceptable accuracy; the calculated value of discrepancy in this case is $\delta_v \approx 0.4$. After spatial filtration, true inhomogeneity $ v^{\text{filtered}}$ and its estimate $ \hat{v}^{\text{filtered}}$ become very close to each other (see figure~\ref{figure_3}d), the discrepancy value is reduced to $\delta_v \approx 0.09$. Relatively high values of $\delta_v$ can be associated with errors in reconstruction of high-frequency spatial components of scatterer function, which, due to the processes of multiple scattering, lead to distortion of spatial frequencies lying in the circle $2k_0$ \cite{label16}. Reconstruction of the considered scatterer in the Born approximation gives unsatisfactory results even after spatial filtering (figure~\ref{figure_3}d); the value of discrepancy in this case is $\delta_v \approx 6.8$.

\begin{figure}[ht!]
	\centerline{\epsfig{file=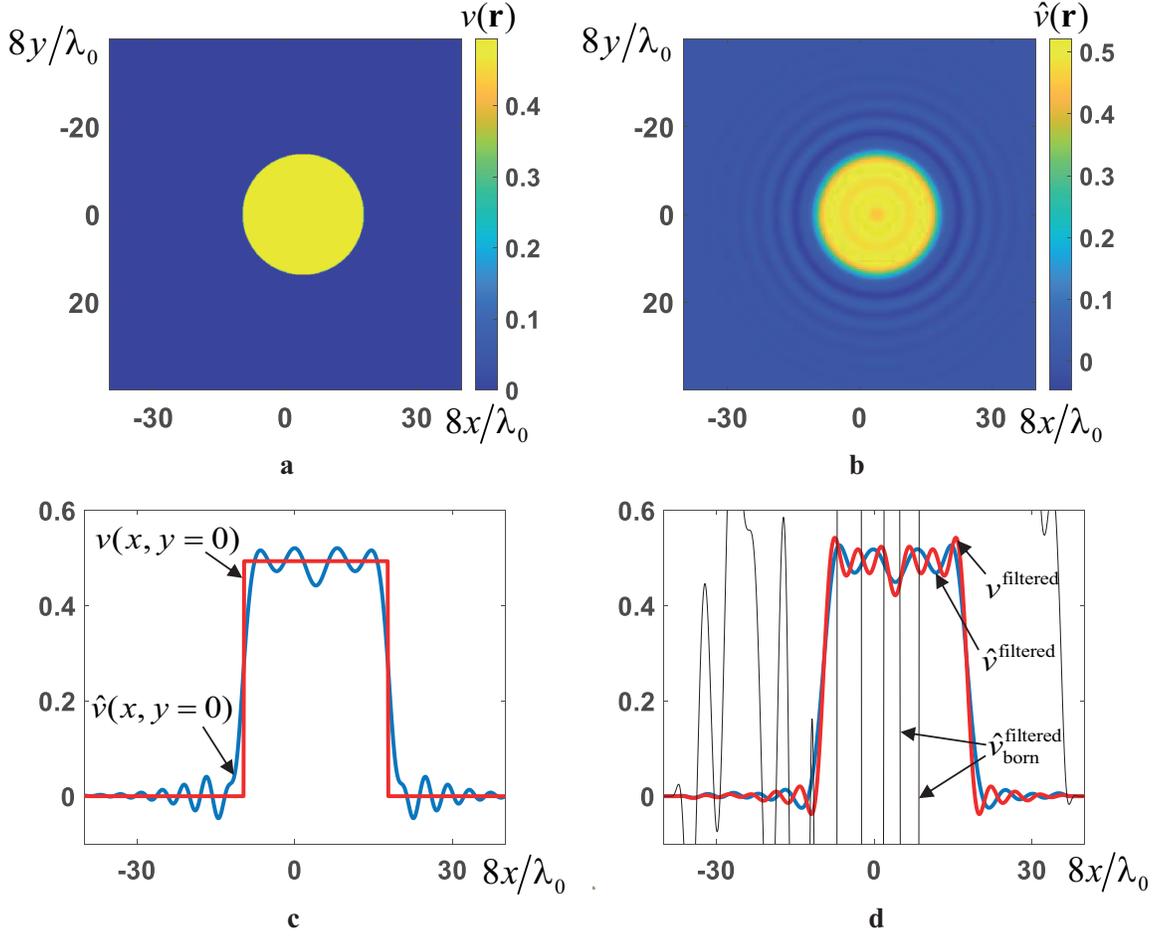, width = 15.7 cm}}
	\caption{Space distribution of true scatterer $v(\mathbf{r})$ (a) with relative contrast of sound speed  $\Delta c \big/ c_0 \approx 1.24$ and corresponding "dielectric" contrast  $\epsilon = 0.2 $, which provides the additional phase shift $\Delta \psi \approx 3.8 \pi$ and norm of scattering data $\| f(\phi, \phi^\prime) \| \approx 18.3 / (3 \pi)$, relative size of inhomogeneity is $ 2 R_0 / \lambda_0 \approx 3.4$;
		\protect\\
		- general view of functional-analytical estimate  $\hat{v}(\mathbf{r})$ (b); 
		\protect\\
		- central cross sections of true inhomogeneity $v(x, y = 0)$ (c, red line) and reconstruction result  $\hat{v}(x, y~=~0)$ (c, blue line); 
		\protect\\
		- central cross sections ($y = 0$) of space filtering results obtained in a cylinder with radius  $2k_0$ (arguments of functions are omitted): filtered true scatterer $ v^{\text{filtered}}$ (d, red line), functional-analytical estimate $ \hat{v}^{\text{filtered}}$ (d, blue line) and Born approximation result $ \hat{v}^{\text{filtered}}_{\text{born}}$ (d, back line).}
	\label{figure_4}
\end{figure}

\begin{figure}[ht!]
	\centerline{\epsfig{file=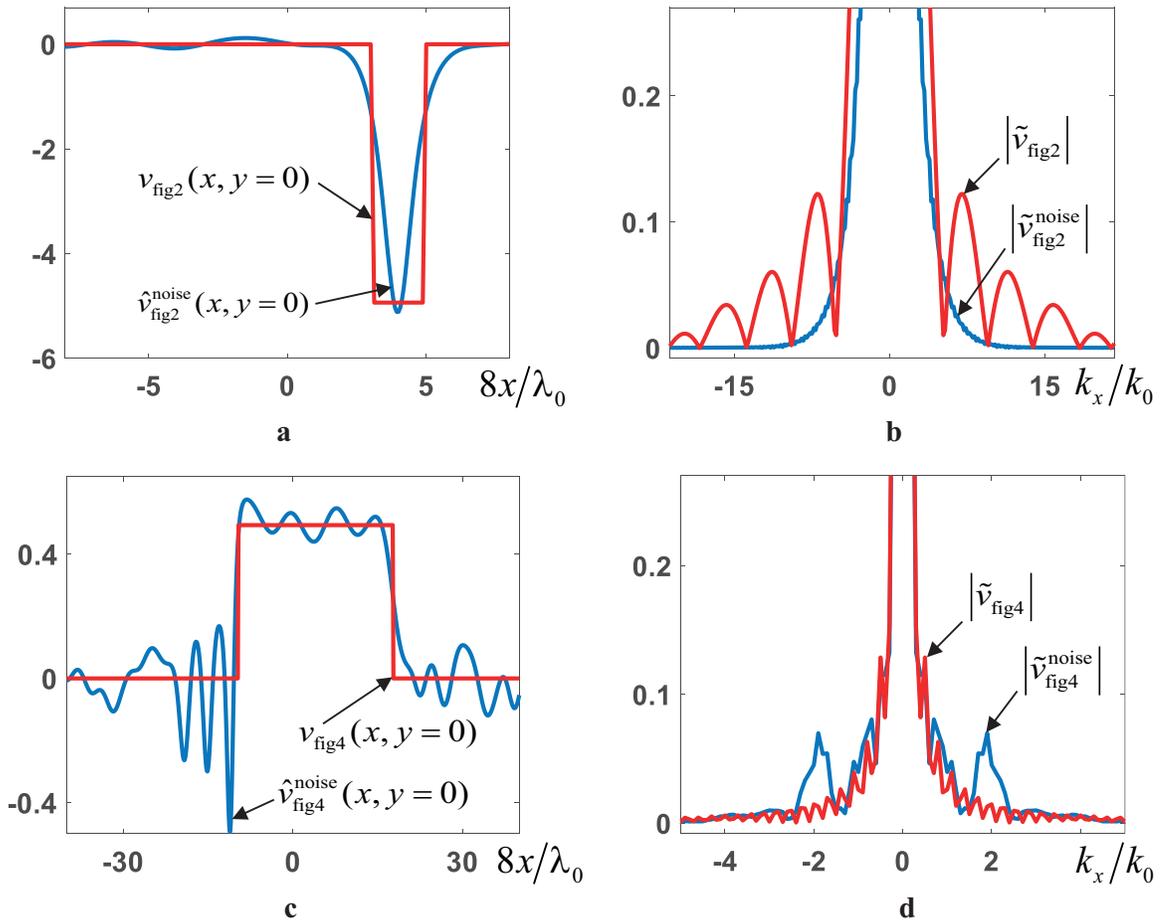, width = 15.7 cm}}
	\caption{Central cross sections of true inhomogeneity $v_\text{fig2}(x, y = 0)$ (a, red line) and reconstruction result  $\hat{v}^{\text{noise}}_\text{fig2}(x, y~=~0)$ (a, blue line) obtained by using noisy data with rms amplitude deviation $\sigma_{\text{ns}} = 0.5 \bar{G}_{\text{sc}}$ (reconstruction of this scatterer with noise-free data is shown in figure~\ref{figure_2});
		\protect\\
		- central cross sections of absolute value of true inhomogeneity space spectrum $ | \tilde v_\text{fig2}(k_x, k_y = 0) |$ normalized on its maximum value (b, red line) and corresponding normalized value of space spectrum modulus of reconstruction result $| \tilde v^{\text{noise}}_{\text{fig2}}(k_x, k_y = 0) |$ (b, blue line), arguments of functions are omitted; 
		\protect\\
		- (c), (d) shows the same as on (a), (b), respectively, but for scatterer function $v_\text{fig4}(\mathbf{r})$ and noisy data with rms amplitude deviation $\sigma_{\text{ns}} = 0.0015 \bar{G}_{\text{sc}}$ (reconstruction of this inhomogeneity with noise-free data is shown in figure~\ref{figure_4}).}
	\label{figure_5}
\end{figure}

The results presented in figure~\ref{figure_3} show capabilities of the functional-analytical algorithm in reconstruction of small contrasts $\epsilon \ll 1$, but correspond to non-physical value of sound speed $c$. In fact, figure~\ref{figure_3} shows reconstruction of the scatterer 
$v(\mathbf{r}) \approx k^2_0$, 
for for limit value 
$\epsilon(\mathbf{r}) \approx 0$, 
where for the considered case 
$k^2_0 = \left( 2 \pi / \lambda_0\right)^2 \approx 0.62$; 
this amplitude of $v(\mathbf{r})$
can be seen in figure~\ref{figure_3}с. In numerical simulation, various contrasts have also been considered in the range $\epsilon \in (0, 1)$ at a fixed value of additional phase shift $\Delta \psi \approx 1.3 \pi$, corresponding to strong scatterer. 
In these tests, instability of reconstruction results, obtained by the functional-analytical algorithm, was not observed; 
the quality of reconstruction for all considered values of $\epsilon$
turned out to be high and close to results shown in figures~\ref{figure_3}b-\ref{figure_3}c. Figure~\ref{figure_4} shows estimates corresponding to more physical values of sound speed in inhomogeneity $c \approx 3.4 \cdot 10^3$~m/s, which gives contrasts $\Delta c \big/ c_0 \approx 1.24$, $\epsilon = 0.2$. This case corresponds, for example, to a metallic object located in water. For the considered problem relative size of inhomogeneity is $ 2 R_0 / \lambda_0 \approx 3.4$ that leads to phase shift $\Delta \psi \approx 3.8 \pi$ and norm of scattering data $\| f(\phi, \phi^\prime) \| \approx 18.3 / (3 \pi)$, i.e. the considered scatterer is even much stronger than it was before in figure~\ref{figure_3}. As it can be seen in figures~\ref{figure_4}a-\ref{figure_4}c, the location of the scatterer function and its amplitude are reconstructed by the functional-analytical algorithm with acceptable accuracy; value of discrepancy is $\delta_v \approx 0.25$ and after spatial filtering reduces to $\delta_v \approx 0.13$. The attempt to reconstruct the considered inhomogeneity in the Born approximation gives unsatisfactory results: only high amplitude oscillations is observed in reconstructed data (figure~\ref{figure_4}d), error after spatial filtering is $\delta_v \approx 6.8$.

Stability of reconstructions shown in figures~\ref{figure_2}-\ref{figure_4} with respect to additive noise in scattering data was also studied.
To make such investigation, a normally distributed random noise 
$
n(\mathbf{y}, \mathbf{x})
$, 
uncorrelated for different emitting and receiving transducers, with zero mean and root-mean-square (rms) amplitude deviation 
$
\sigma_{\text{ns}} = \alpha \bar{G}_{\text{sc}}
$,
was introduced into the scattered fields
$
G_{\text{sc}}(\mathbf{y}, \mathbf{x}) \equiv G(\mathbf{y}, \mathbf{x}) - G_0(\mathbf{y}, \mathbf{x})
$
separately for real and imaginary parts; 
here $\mathbf{x}$, $\mathbf{y}$ are radius vectors of source and receiver, respectively, and 
the rms value $\bar{G}_{\text{sc}}$ is defined as 
$
\bar{G}_{\text{sc}} \equiv 
\sqrt{
	\int\limits_{S} d \mathbf{x} \int\limits_{S} d \mathbf{y} \ 
	\bigl| G_{\text{sc}}(\mathbf{y}, \mathbf{x}) \bigr|^2 } 
\bigg/ 
\sqrt{
	\int\limits_{S} d \mathbf{x} \int\limits_{S} d \mathbf{y} } 
$.

Numerical investigation implemented with different noise levels has revealed that reconstructions shown in figures~\ref{figure_2},~\ref{figure_3} are very stable even for $ \alpha \approx 0.5 $. On the other hand, the reconstruction shown in figure~\ref{figure_4} turned out to be unstable even for very small noise level with $ \alpha \approx 0.002 $. 
As an example, figure~\ref{figure_5} shows reconstruction results obtained from noisy scattering data for inhomogeneities shown in figures~\ref{figure_2},~\ref{figure_4}; these inhomogeneities will be referred bellow as $v_\text{fig2}(\mathbf{r})$ and $v_\text{fig4}(\mathbf{r})$, respectively. 
As it can be seen in figure~\ref{figure_5}a, reconstruction $\hat{v}^{\text{noise}}_\text{fig2}(\mathbf{r})$ 
obtained with noise level $\alpha = 0.5$ is almost the same as in figure~\ref{figure_2}; location of the scatterer function and its amplitude are reconstructed with acceptable accuracy, value of discrepancy after spatial filtering is $\delta_v \approx 0.18$. At the same time, estimate 
$\hat{v}^{\text{noise}}_\text{fig4}(\mathbf{r})$ 
reconstructed with very small level of noise $\alpha = 0.0015$ has strong distortions (see figure~\ref{figure_5}c) in comparison with figure~\ref{figure_4}c; discrepancy after spatial filtering is $\delta_v \approx 0.26$.
The considered result can be explained in terms of spatial spectrum localization of so-called secondary sources, which are equal to $ p(\mathbf{r}) v(\mathbf{r})$ \cite{label1, label16}. 
It is known \cite{label1, label16} that due to multiple scattering of incident wave on strong scatterers, spatial spectrum of secondary sources can expands even beyond the circle of radius $ 2k_0$. 
It should be noted that scatterer shown in figure~\ref{figure_4} is much stronger ($\Delta \psi \approx 3.8 \pi$) than inhomogeneities shown in figures~\ref{figure_2} and~\ref{figure_3} ($\Delta \psi \lesssim 1.3 \pi$), i.e. processes of multiple scatterings are more pronounced for $ v_\text{fig4}(\mathbf{r})$.
For small-scale inhomogeneity $v_\text{fig2}(\mathbf{r})$ its spatial spectrum $\tilde{v}_\text{fig2}(\mathbf{k})$ (i.e. Fourier transform of $v_\text{fig2}(\mathbf{r})$)
is wide; as a result, high spatial frequencies, for example, near $ 2k_0 $, are not suppressed by noise (see figure~\ref{figure_5}b).
Since high spatial frequencies determine reconstruction of small-scale inhomogeneities, the scatterer $v_\text{fig2}(\mathbf{r})$ is reconstructed with good noise immunity (see figure~\ref{figure_5}a); the similar results were obtained for scatterer shown in figure~\ref{figure_3}.
The inhomogeneity $v_\text{fig4}(\mathbf{r})$ has narrow spatial spectrum, so noise dominates at high frequencies (see figure~\ref{figure_5}d).
Strong noise in high-frequency components of spatial spectrum, in combination with multiple scattering processes, leads to strong distortions of reconstruction (see figure~\ref{figure_5}c).
As a result, the relatively weak scatterers shown in figures~\ref{figure_2} and~\ref{figure_3} with wide spatial spectrum turned out to have better noise immunity compared to the strong scatterer shown in figure~\ref{figure_4}, which has narrow spatial spectrum.
It should be noted that the functional-analytical algorithm allows to reconstruct inhomogeneities by using scattering data measured at many frequencies \cite{label17} that can be implemented in practice by using multifrequency or impulse sounding regime. This is a very important way to improve results of reconstruction from noisy data even in the case of strong scatterers with complicated space distribution of its parameters \cite{label16, label18}.

\renewcommand{\thesection}{\large 5}
\section{\large Conclusions}
The results obtained in the present work show that using the considered functional-analytical algorithm \cite{label14,label15} it is possible to reconstruct acoustic inhomogeneities with high equivalent "dielectric" contrast $\epsilon \gg 1$. This, in some sense, demonstrates the advantages of functional-analytical approach in comparison with least squares methods of inverse problem solution. Moreover, the functional algorithm also perfectly reconstructs very small "dielectric" contrast  $\epsilon \ll 1$, recovering of which can be difficult for other approaches due to small backscattering for such weak contrasts.

It should be noted that to compare reconstruction results of different scatterers obtained by different methods it makes sense to focus not on contrasts but on the so-called scatterer strength \cite{label1}, which characterizes how strongly the incident field is distorted during scattering. If scatterer is described only by sound speed perturbation, then an additional phase shift $\Delta \psi$ is used as a quantitative characteristic of scatterer strength \cite{label1}; in the case of inhomogeneity with cylindrical shape and radius $R_0$, the simple relation can be used $\Delta \psi = 2 R_0 \left[ k_0 - k\right]$, where $k_0$ and $k$ are wave numbers in background medium and in inhomogeneity, respectively. If scatterer is strong, then values of $| \Delta \psi |$ is of order or greater than $\pi$ \cite{label1}. Moreover, a scatterer with $\Delta \psi < 0$ appears to be more difficult to reconstruct, due to its focusing properties (focusing scatterer), than an inhomogeneity with $\Delta \psi > 0$ (defocusing scatterer) \cite{label33}. An alternative estimate of the scatterer strength is the norm of scattering data $\| f \|$. In reconstructions of the present work the greatest value of scattering data norm is $\| f \| \approx 18.3 / (3 \pi)$. However, this is not a limit for the functional-analytical algorithm. For example, in \cite{label29}, reconstruction results with high quality were obtained from scattering data with norm $\| f \| \approx 19.3 / (3 \pi)$. Even stronger scatterers can be reconstructed if spatial spectrum of secondary sources $v(\mathbf{r}) p(\mathbf{r})$ is localized inside the circle of radius $2k_0$ \cite{label1,label16} and initial scattering data are obtained with a large number of emission and reception points. In this case, the functional-analytical algorithm gives reconstruction results with accuracy depending mainly on accuracy of initial scattering data, rather than on inhomogeneity contrasts \cite{label16}. 


\section{\large Acknowledgement}

The author thanks M.V. Klibanov for the suggestion to carry out numerical simulation, results of which are presented in this work. The author also thanks R.G. Novikov for valuable discussions.
\\
\\
The reported study was funded by RFBR and CNRS, project number 20-51-15004.


\end{document}